\documentclass[11pt]{article}
\usepackage{enumerate}
\usepackage{tikz}
\usepackage{amssymb,a4wide,latexsym,makeidx,epsfig,fleqn}
\usepackage{amsthm}
\usepackage{amsmath}
\usepackage{enumerate}
\usepackage{array}
\newpage
\newtheorem{theorem}{Theorem}[section]
\newtheorem{remark}[theorem]{Remark}

\newtheorem{definition}[theorem]{Definition}
\newtheorem{lemma}[theorem]{Lemma}

\newtheorem{corollary}[theorem]{Corollary}

\def\Tr{\mathop{\rm Tr }\nolimits}
\def\tr{\mathop{\rm tr }\nolimits}

\def\Spec{\mathop{\rm Spec }\nolimits}

\begin{document}
\textwidth 150mm \textheight 225mm
\title{The characteristic polynomials of $r$-uniform hypercycles with length $l$
\thanks{Supported by the National Natural Science Foundation of China (Nos. 12271439 and 12301452) and the Natural Science Foundation of Shanxi Province (No. 202303021212025).}}
\author{{Bo Dong$^{a,b}$, Cunxiang Duan$^{c}$, Ligong Wang$^{a,b,}$\footnote{Corresponding author.}}\\
{\small $^{a}$School of Mathematics and Statistics, Northwestern
Polytechnical University,}\\ {\small Xi'an, Shaanxi 710129, P.R. China}\\
{\small $^{b}$Xi'an-Budapest Joint Research Center for Combinatorics, Northwestern Polytechnical University,}\\ {\small Xi'an, Shaanxi 710129, P.R. China}\\
{\small $^{c}$School of Mathematics, Taiyuan University of Technology,}\\ {\small Taiyuan, Shanxi 030024, P.R. China}\\
\\{\small E-mail: dongbo98@163.com; cxduanmath@163.com; lgwangmath@163.com}}
\date{}
\maketitle
\begin{center}
\begin{minipage}{120mm}
\vskip 0.3cm
\begin{center}
{\small {\bf Abstract}}
\end{center}
{\small Let $C_{l}$ be a cycle with length $l.$ The $r$-uniform hypercycle with length $l$ is obtained by adding $r-2$ new vertices in every edge of $C_{l},$ denoted by $C_l^{(r)}$. In this paper, we deduce some higher-order traces for the adjacent tensor of $C_l^{(r)}$ by BEST Theorem. Then we obtain higher-order spectral moments according to the relationship between eigenvalues of power hypergraphs and eigenvalues of signed graphs. Finally, the general expression of the characteristic polynomials of $C_l^{(r)}$ is given. Furthermore, by using this general expression,  we present the characteristic polynomials of $C_5^{(r)}$ and $C_6^{(r)}$ as examples.
}
\end{minipage}
\end{center}

\section{Introduction }
\label{sec:ch6-introduction}

Denote the $n$-element set $\{1,2,\ldots,n\}$ by $[n]$. Let $H=(V(H),E(H))$ be a hypergraph with vertex set $V(H)=[n]$ and hyperedge (abbreviated as edge) set $E(H)$. Specially, if every edge $e\in E(H)$ has the same cardinality $r$ (i.e. $|e|=r$), then we call $H$ is a $r$-uniform hypergraph ($r$-graph for short). Note that a $2$-graph is actually a graph. For a graph $G=(V(G),E(G))$, denote by $G^{(r)}$ the $r$-uniform hypergraph obtained by adding $r-2$ new vertices to each edge of $G$, which is called the $r$-power hypergraph of $G$. Thus $|V(G^{(r)})|=|V(G)|+(r-2)|E(G)|$.

An $r$-order $n$ dimensional complex tensor $\mathcal{T}=(t_{i_{1}i_{2}\ldots i_{r}})$ is a multidimensional array with $n^r$ complex entries, where $i_j\in[n]$, $j=1,2,\ldots,r$. In 2012, Cooper and Dulte \cite{CoDu} introduced the concept of adjacency tensor $\mathcal{A}(H)$ for an $r$-uniform hypergraph $H$. Spectra of hypergraphs are based on the spectral theory of tensors. Throughout this paper, the characteristic polynomial (resp., the spectrum and the eigenvalues) of an uniform hypergraph refers to that of the adjacent tensor of the hypergraph. Recently, the spectral theory of hypergraphs has emerged as a widely researched and deeply analyzed attracting subject in spectral graph theory. A lot of results 
\cite{Bret,Qibook,SSZ} provided the  foundation and tools for subsequent research on spectra of hypergraphs.

The study of characteristic polynomials of graphs is a very classical problem in graph theory. However, that of hypergraphs is much more novel and challenging since the spectra of hypergraphs based tensors is more difficult to calculate. Up to now, the range of hypergraphs for which explicit expressions of characteristic polynomials have been explored remains limited. Much of the research has focused on hypertrees, and the Poisson Formula has played an important role. Initially, Cooper and Dutle \cite{CoDu} gave the characteristic polynomial of a single hyperedge by the characteristic equation in 2012. In 2015, Cooper and Dutle \cite{CoD} obtained the characteristic polynomials of 3-uniform hyperstars by using the Poisson Formula. Further, in 2018, Bao et al. \cite{BFWZ} improved this result to $r$-uniform hyperstars. In 2019, Chen and Bu \cite{CB}, derived a reduced formula by the Poisson Formula to calculate  the characteristic polynomials of hypergraphs with pendant edges. By using this reduced formula, they presented the characteristic polynomials of $r$-uniform hyperpaths. In 2023, Duan et al. \cite{D2023} obtained the characteristic polynomials of uniform double hyperstars. In 2024, Chen et al. \cite{CDB2} determined the entire spectrum of power hypergraphs interms of parity-closed walks. In 2024, Li et al. \cite{hypertree} revealed the equality between the characteristic and matching polynomials of hypertrees, and gave a general expression for the characteristic polynomials of hypertrees.

Besides acyclic hypergraphs (i.e. hypertrees), the characteristic polynomials have been obtained only for some well-structured hypergraphs. In 2021, Zheng \cite{Z} gave the characteristic polynomials of complete 3-uniform hypergraphs. In 2023 and 2024, Duan et al. \cite{D2023,D2024} obtained the characteristic polynomials of $r$-uniform hypercycles with lengths 3 and 4, respectively. In this process, the method of studying hypergraph spectra through higher-order trace formulas for tensors \cite{MS,SQH}  has played a very important role.

Another approach to studying the characteristic polynomial of hypergraphs is through its coefficients. In this regard, there is the generalized Harary-Sachs theorem \cite{CC} for uniform hypergraphs.

For other results about spectra of hypergraphs, please refer to \cite{CBZ,CoDu,DW,GCH,HHLQ,KLQY,ZSWB}

In this paper, we mainly study the characteristic polynomial of an $r$-uniform hypercycle with length $l$. In Section 2, we focus on introducing concepts related to the spectra, higher-order traces of tensors (resp. uniform hypergraphs), and briefly outline the derivation of this article. In Section 3, we give some high-order traces of $\mathcal{A}(C_{l}^{(r)})$. In Section 4, according to the relationship between high-order traces and spectral moments, we solve for the unknown parameters in the characteristic polynomial of $C_{l}^{(r)}$, thereby obtaining the general formula of characteristic polynomial of $C_{l}^{(r)}$. Finally, according to the general formula, we present the characteristic polynomials of $C_{5}^{(r)}$ and $C_{6}^{(r)}$ as examples.


\section{Preliminaries}
\label{sec:ch-sufficient}

Denote by $\mathbb{C}$ (resp. $\mathbb{R}$) the field of complex (resp. real) numbers. Let $\mathbb{C}^{[r,n]}$ be the set of all $r$-order $n$ dimensional complex tensors. In 2013, Shao \cite{Shao} introduced a general product of tensors as follows.

\noindent\begin{definition}\label{de:product} (\cite{Shao})
Let $\mathcal{A}=(a_{i_1i_2\cdots i_r})\in \mathbb{C}^{[r,n]}$ and $\mathcal{B}=(b_{i_1i_2\cdots i_t})\in \mathbb{C}^{[t,n]}$, where $r\geq 2$ and $t\geq 1$. Define the product $\mathcal{A}\cdot \mathcal{B}=\mathcal{A}\mathcal{B}$ to be a complex tensor of order $(r-1)(t-1)+1$ and dimension $n$ as follows
\begin{equation*}
    c_{i\alpha_1\alpha_2\cdots \alpha_{r-1}}=\sum_{i_2,i_3,\ldots,i_r=1}^{n}a_{ii_2i_3\cdots i_r}b_{i_2\alpha_1}b_{i_3\alpha_2}\cdots b_{i_r\alpha_{r-1}}~(i\in[n],~ \alpha_1,\alpha_2,\ldots,\alpha_r\in[n]^{t-1}).
\end{equation*}
\end{definition}
Specially, when $t=1$, $\mathcal{B}=\boldsymbol{x}\in\mathbb{C}^n$ is an $n$ dimensional complex vector. Then $\mathcal{A}\mathcal{B}=\mathcal{A}\boldsymbol{x}$ is an $n$ dimensional vector with the following $i$-th component
\begin{equation*}
    (\mathcal{A}\boldsymbol{x})_i=\sum_{i_2,i_3,\ldots,i_r=1}^{n}a_{ii_2i_3\cdots i_r}x_{i_2}x_{i_3}\cdots x_{i_r}.
\end{equation*}

In 2005, Qi \cite{Qi} and Lim \cite{Lim} proposed the eigenvalues of tensors independently. Let $\mathcal{I}=(\delta_{i_1i_2\cdots i_r})\in \mathbb{C}^{[r,n]}$, where $\delta_{i_1i_2\cdots i_r}=1$ if and only if $i_1=i_2=\cdots =i_r$, otherwise $\delta_{i_1i_2\cdots i_r}=0$.

\noindent\begin{definition}\label{de:eigenvalue} (\cite{Qi,Lim})
    Let $\mathcal{A}\in \mathbb{C}^{[r,n]}$. A pair $(\lambda,\boldsymbol{x})\in \mathbb{C}\times(\mathbb{C}^{n}\setminus \{0\})$ is called an eigenvalue and eigenvector corresponding to $\lambda$ of $\mathcal{A}$ if they satisfy the following homogeneous polynomial equation (i.e. the so called characteristic equation)
    \begin{equation*}
        \mathcal{A}\boldsymbol{x}=\lambda\mathcal{I}\boldsymbol{x}.
    \end{equation*}
\end{definition}

Here $(\mathcal{A}\boldsymbol{x})_i$ is exactly $(\mathcal{A}\boldsymbol{x}^{r-1})_i$ in \cite{Qi}. And $(\mathcal{I}\boldsymbol{x})_i=x_i^{r-1}$ is identified with $\boldsymbol{x}^{[r-1]}$ in \cite{Qi}.

Initially, Qi \cite{Qi} introduced the concept of symmetric hyperdeterminants (determinant for shorts) to investigate the eigenvalues of symmetric tensors. The concept is based on resultants. And for more information about resultants, see the paper \cite{MS} of Morozov and Shakirov, in which higher-order traces were introduced to give an expression of resultant. For a tensor $\mathcal{T}=(t_{i_1i_2\cdots i_r})\in \mathbb{C}^{[r,n]}$, denote by $Det(\mathcal{T})$ (resp. $Res(\mathcal{T}\boldsymbol{x})$) its determinant (resp. the resultant of the system of polynomial equations $\mathcal{T}\boldsymbol{x}=0$). Then according to \cite{HHLQ} by Hu et al.,
\begin{equation*}
    Det(\mathcal{T})=Res(\mathcal{T}\boldsymbol{x}).
\end{equation*}

\noindent\begin{definition}\label{de:defofcharpoly}
The characteristic polynomial of $\mathcal{T}$, denoted by $\phi(\mathcal{T};\lambda)$, is a monic polynomial in $\lambda$ of degree $n(r-1)^{n-1}$, whose roots are excatlly all eigenvalues of $\mathcal{T}$, specifically
\begin{equation*}
    \phi(\mathcal{T};\lambda)=Det(\lambda\mathcal{I}-\mathcal{T})=Res((\lambda\mathcal{I}-\mathcal{T})\boldsymbol{x}).
\end{equation*}
\end{definition}

Let $\mathcal{T}=(t_{i_1i_2\cdots i_r})\in \mathbb{C}^{[r,n]}$ and $A=(a_{ij})$ be an auxiliary matrix. Define corresponding differential operators as follows
\begin{equation*}
    \hat{g}_i=\sum_{i_2=1}^n\sum_{i_3=1}^n\cdots\sum_{i_r=1}^n t_{ii_2i_3\cdots i_r}\frac{\partial}{\partial a_{ii_2}}\frac{\partial}{\partial a_{ii_3}}\cdots\frac{\partial}{\partial a_{ii_r}}, ~i\in[n].
\end{equation*}
\noindent\begin{definition}\label{trace1}(\cite{MS})
For a tensor $\mathcal{T}=(t_{i_1i_2\cdots i_r})\in \mathbb{C}^{[r,n]}$, its $d$-th order trace is defined as
\begin{equation*}
        \Tr_{d}(\mathcal{T})=(r-1)^{n-1}\bigg [\sum_{d_{1}+d_{2}+\cdots+d_{n}=d}\prod_{i=1}^{n}\frac{(\hat{g}_i)^{d_i}}{(d_{i}(k-1))!} \bigg ]\tr (A^{d(r-1)}),
\end{equation*}
where every entries of auxiliary matrix $A=(x_{ij})$ will be eliminated by the differential operators and $\tr$ refers to the trace of a matrix.
\end{definition}
Then by Morozov and Shakirov \cite{MS}, the characteristic polynomial can be expressed as follows
\begin{equation*}
    \phi(\mathcal{T};\lambda)=\sum_{j=0}^{N}P_{j}(-\frac{\Tr_{1}(\mathcal{A})}{1}, -\frac{\Tr_{2}(\mathcal{A})}{2}, \ldots, -\frac{\Tr_{j}(\mathcal{A})}{j})\lambda^{N-j},
\end{equation*}
where $N=n(r-1)^{n-1}$ and
\begin{equation*}
    P_0(t_0)=1;\qquad P_j(t_1,t_2,\ldots,t_j)=\sum_{m=1}^{j}\sum_{h_1+h_2+\cdots+h_m=j}\frac{t_{h_1}t_{h_2}\cdots t_{h_i}}{m!}~\text{for}~ j\geq 1.
\end{equation*}

\noindent\begin{theorem}\label{trace=moment}(\cite{HHLQ})
    For a tensor $\mathcal{T}=(t_{i_1i_2\cdots i_r})\in \mathbb{C}^{[r,n]}$, denote by $Spec(\mathcal{T})$ the multiset of eigenvalues of $\mathcal{T}$ (i.e. the spectrum of $\mathcal{T}$), then
    \begin{equation*}
        \Tr_d(\mathcal{T})=\sum_{\lambda\in Spec(\mathcal{T})}\lambda^{d}\quad\text{for~}1\leq d\leq N=n(r-1)^{n-1}.
    \end{equation*}
    The right hand of the equation is called the $d$-th spectral moment of $\mathcal{T}$.
\end{theorem}

Based on the hight-order trace gave by Morozov and Shakirov \cite{MS}, Shao et al. \cite{SQH} presented a formula for higher-order traces that does not rely on the auxiliary matrix or differential operators, but is based on multi-digraphs (digraphs with multi-arcs).

For a positive integer $d$, let
\begin{equation*}
    \mathcal{F}_d=\left\{(i_1\alpha_1,i_2\alpha_2,\ldots,i_d\alpha_d):1\leq i_1\leq i_2\alpha_2\leq\cdots\leq i_d\leq n,~\alpha_1,\alpha_2,\ldots,\alpha_d\in [n]^{r-1}\right\}.
\end{equation*}
For $F=(i_1\alpha_1,i_2\alpha_2,\ldots,i_d\alpha_d)\in \mathcal{F}_d$ and $\mathcal{T}=(t_{i_1i_2\cdots i_r})$, define a function $\pi_{F}(\mathcal{T})=\prod_{j=1}^{d}t_{i_j\alpha_{j}}$. We call $F$ is $r$-valent if every index appearing as some entry of $F$ appears $r$-multiple times. For a digraph $D=(V,A)$, the in-degree (resp. out-degree) of $v\in V$ is the number of arcs incident to (resp. from) $v$, denoted by $d^{-}(v)$ (resp., $d^{+}(v)$).

\noindent\begin{definition}\label{de:trace2} (\cite{SQH})
Let $F=(i_{1}\alpha_{1},i_2\alpha_2,\ldots,i_{d}\alpha_{d})\in \mathcal{F}_{d}$, where $\alpha_s=v_{1}^{(s)}v_2^{(s)}\cdots v_{r-1}^{(s)}\in [n]^{r-1}$ for $s\in[d]$. Define that
\begin{enumerate}
    \item[(1)] $E_s(F)=\left\{(i_s,v_{t}^{(s)}):t\in[r-1]\right\}$ and $E(F)=\bigcup_{s=1}^d E_s$ are both arc multi-set. $D(F)=D=(V(F),E(F))$ is the multi-digraph corresponding to arc multi-set $E(F)$, where $V(F)$ is the set of indices appearing in $F$.
    \item[(2)] $b(F)=\prod\limits_{a\in E(F)} m(a)!$, $c(F)=\prod\limits_{v\in V(F)}d^{+}(v)!,$ where $m(a)$ is the multiplicity of the arc $a.$
    \item[(3)] $W(F)$ is the set of all Eulerian closed walks in digraph $D$.
\end{enumerate}
Note that multi-arcs of $W(F)$ are not distinguished and $W(F)=\emptyset$ if $F$ is not $r$-valent. Then the $d$-th order trace of $\mathcal{T}\in\mathbb{C}^{[r,n]}$ can be defined equivalently as
\begin{equation*}
    \Tr_{d}(\mathcal{T})=(r-1)^{n-1}\sum_{F\in \mathcal{F}'_{d}}\frac{b(F)}{c(F)}\pi_{F}(\mathcal{T})\lvert W(F)\rvert,
\end{equation*}
where $\mathcal{F}'_{d}=\{F\in \mathcal{F}_{d} : F$ is $r$-valent$\}.$
\end{definition}

\noindent\begin{definition}\label{de:c1} (\cite{CoDu})
For a $r$-uniform hypergraph $H = (V(H), E(H))$ with $V(H)=[n]$, its adjacency tensor, denoted by $\mathcal{A}(H)=(a_{i_{1}i_{2}\ldots i_{r}})$, is defined as follows
\begin{equation*}
a_{i_{1}i_{2}\ldots i_{r}}=\left\{
\begin{array}{cl}
\frac{1}{(r-1)!}& \mbox {if}   ~\{i_{1},i_{2},\ldots, i_{r}\} \in E,
\\
0& \mbox {otherwise}.
\end{array}
\right.
\end{equation*}
\end{definition}

For a graph $G$, let $\pi: E\rightarrow \{+1,-1\}$ be an edge sign function, then the pair $(G,\pi)$ is called a signed graph, denoted by $G_{\pi}.$ The edge $e\in E(G)$ is called negative edge if $\pi(e)=-1$. An (induced) subgraph of $G_{\pi}$ is called a signed (induced) subgraph of $G$. For two vertices $i,j\in V(G)$, denote by $i\sim j$ (resp. $i\nsim j$) that $i$ and $j$ are adjacent in $G$ (resp. $i$ and $j$ are not adjacent). Then adjacent matrix $A(G_{\pi})=(a_{ij})$ is defined as
\begin{equation*}
    a_{ij}=\left\{
    \begin{array}{cl}
        \pi(i,j) & \mbox{if~} i\sim j, \\
         0 & \mbox{if~} i\nsim j.
    \end{array}
    \right.
\end{equation*}In 2023, Chen et al. \cite{CDB} revealed the relationship between eigenvalues of signed (induced) subgraph and the $r$-power $G^{(r)}$ of $G$.

\noindent\begin{lemma}\label{le:4-5}(\cite{CDB}) $\lambda\in \mathbb{C}$ is an eigenvalue of the $r$-power $G^{(r)}$ of a graph $G$ if and only if
\begin{enumerate}
    \item[(1)] for $r = 3,$ $\beta$ is an eigenvalue of some signed induced subgraph of $G$ and $\beta^{2}= \lambda^{r}$,
    \item[(2)] for $r\geq 4,$ $\beta$ is an eigenvalue of some signed subgraph of $G$ and $\beta^{2}= \lambda^{r}.$
\end{enumerate}
\end{lemma}

\noindent\begin{lemma}\label{le:r-symmetry}(\cite{SQH})
Let $H$ be an $r$-uniform hypergraph, and $\phi(H;\lambda)=\sum_{j=0}^{N}a_j\lambda^{N-j}$ ($N=n(r-1)^{n-1}$) be its characteristic polynomial. Then the following three condition holds and are equivalent
\begin{enumerate}
    \item[(1)] The spectrum of $H$ is $r$-symmetric.
    \item[(2)] If $d$ is not a multiple of $r$, then $a_d=0$, namely, there exists some integer $t$ and polynomial $f$ such that $\phi(H;\lambda)=\lambda^t f(\lambda^r)$.
    \item[(3)] if $d$ is not a multiple of $r$, then $\Tr_d(\mathcal{A}(H))=0$.
\end{enumerate}
\end{lemma}

Now we focus on all the $r$-th powers of eigenvalues of $r$-uniform hypercycle $C_{l}^{(r)}$, i.e., all squared eigenvalues of signed subgraph of $C_l$.

\noindent\begin{definition}\label{de:+-cycles}
A signed cycle is said to be positive if it contains an even number of negative edges,
otherwise the cycle is called negative. A signed graph is said to be balanced if all its cycles are positive, otherwise it is unbalanced.
\end{definition}

In the spectral theory of signed graphs, an important concept is switching equivalence. Two switching equivalent signed graphs share identical spectra and characteristic polynomials. For more details of the spectral properties of signed graphs and switching equivalence, refer to references \cite{ABD,ZT}

\noindent\begin{lemma}(\cite{ABD})\label{le:ABD}
If $G$ is a tree, then all signed graphs on $G$ are switching equivalent, i.e., the spectra of these signed graphs are identical with the spectrum of $G$; If $G$ is a cycle $C_l$, then there are exactly two switching equivalent classes for signed graphs on $G$: the positive cycles and the negative cycles.
\end{lemma}

\noindent\begin{lemma}\cite{CvRS}\label{le:spectra}
Let $C_n$ (resp. $P_n$) be the $n$-vertex cycle (resp. path), then
\begin{align*}
    \Spec(C_n) &= \{2\cos(\frac{2k\pi}{n}):k\in[n]\},\\
    \Spec(P_n) &= \{2\cos(\frac{k\pi}{n+1}):k\in[n]\}.
\end{align*}
\end{lemma}

\noindent\begin{lemma}\label{le:r-th}
Let $C_{l}^{(r)}$ be an $r(\geq 4)$-uniform hypercycle with length $l$, then its spectrum $\Spec(C_{l}^{(r)})$ satisfies
\begin{equation*}
    \Big\{\lambda^r:\lambda\in \Spec(C_{l}^{(r)})\Big\}=\{4\}\bigcup \Big\{(2\cos(\frac{k\pi}{j+1}))^2:k\in[j],j\in[l]\Big\}.
\end{equation*}
\end{lemma}

\noindent\textbf{Proof.} According to Lemmas \ref{le:4-5}, \ref{le:ABD} and \ref{le:spectra}, the squared eigenvalues of all signed subpath are listed in Table \ref{tab:kpowereignvalues}.

According to Lemma \ref{le:ABD}, signed cycles with length $l$ are divided into two classes, positive and negative signed cycles with length $l$. The spectra of positive signed $C_l$'s are identical with the spectrum of $C_l$, shown in Lemma \ref{le:spectra}. After some calculations, the spectra of negative signed $C_l$'s are the same as
\begin{equation*}
    \Big\{2\cos(\frac{(2k-1)\pi}{l}):k\in[l]\Big\},
\end{equation*}
which is contained in $\Spec(P_{l-1})$. Except for $(2\cos(\frac{2l\pi}{l}))^2=4$, a squared eigenvalue of the positive signed $C_l$, all other squared eigenvalues of signed $C_l$ are included in Table \ref{tab:kpowereignvalues}. $\hfill$$\square$

\begin{table}[h]
	\centering
	\begin{tabular}{lllcl}
		\cline{1-2}
		$P_1$: & \multicolumn{1}{|l}{$(2\cos(\frac{\pi}{2}))^2=0$} & \multicolumn{1}{|l}{}                              & \multicolumn{1}{l}{}        &                                                     \\ \cline{1-3}
		$P_2$: & \multicolumn{1}{|l}{$(2\cos(\frac{\pi}{3}))^2$}   & \multicolumn{1}{|l}{$(2\cos(\frac{2\pi}{3}))^2$}   & \multicolumn{1}{|l}{}       &                                              \\ \cline{1-4}
        $P_3$: & \multicolumn{1}{|l}{$(2\cos(\frac{\pi}{4}))^2$}   & \multicolumn{1}{|l}{$(2\cos(\frac{2\pi}{4}))^2$}   &  \multicolumn{1}{|l}{$(2\cos(\frac{3\pi}{4}))^2$}  & \multicolumn{1}{|l}{}        \\

		\multicolumn{5}{c}{$\vdots$}                                                                                                                                                                          \\ \hline
		$P_l$: & \multicolumn{1}{|l}{$(2\cos(\frac{\pi}{l+1}))^2$} & \multicolumn{1}{|l}{$(2\cos(\frac{2\pi}{l+1}))^2$} & \multicolumn{1}{|c}{$\cdots$} & \multicolumn{1}{|l|}{$(2\cos(\frac{l\pi}{l+1}))^2$} \\ \hline
	\end{tabular}
	\caption{All squared eigenvalues of signed paths}
	\label{tab:kpowereignvalues}
\end{table}

\noindent\begin{remark}\label{re:setr=3}
By Lemma \ref{le:4-5}, when $r=3$, the spectrum $\Spec(C_{l}^{(3)})$ satisfies
\begin{equation*}
    \Big\{\lambda^3:\lambda\in \Spec(C_{l}^{(3)})\Big\}=\{4\}\bigcup \Big\{(2\cos(\frac{k\pi}{j+1}))^2:k\in[j],j\in[l-1]\Big\}.
\end{equation*}
\end{remark}

Throughout this paper, let $\lambda_{1,1}=2$ and $\lambda_{j,k}=2\cos(\frac{k\pi}{j+1})$ for $2\leq j\leq l$ and $1\leq k\leq j$. For any $n$-degree monic polynomial $f(\lambda)=\prod_{i=1}^{n}(\lambda-\lambda_i)$, a transformation $f(\lambda)\mapsto \tilde{f}(\lambda)$ on monic polynomials is defined such that $\tilde{f}(\lambda)=\prod_{i=1}^{n}(\lambda^r-\lambda_i^2)$.
\noindent\begin{lemma}
Let $C_{l}^{(r)}$ be an $r(\geq 4)$-uniform hypercycle with length $l$. Then $\phi(C_{l}^{(r)};\lambda)$ has the following form
\begin{equation}\label{eq:1}
    \phi(C_{l}^{(r)};\lambda)=\lambda^{m_0}(\lambda^r-\lambda_{1,1}^2)^{m_1}\prod_{j=2}^{l}[\tilde{\phi}(P_j;\lambda)]^{m_j},
\end{equation}
where $\tilde{\phi}(P_j;\lambda)=\prod_{k=1}^{j}(\lambda^r-\lambda_{j,k}^2)$ is the transformed polynomial of $\phi(P_{j};\lambda)$.
\end{lemma}

\noindent\textbf{Proof.} According to Lemmas \ref{le:4-5} and \ref{le:r-symmetry}, $\phi(C_{l}^{(r)};\lambda)$ has the following form
\begin{equation*}
    \phi(C_{l}^{(r)};\lambda)=\lambda^{m_0}\prod_{j=1}^{l}\prod_{k=1}^{j}(\lambda^r-\lambda_{j,k}^2)^{m(j,k)}.
\end{equation*}

According to Theorem \ref{trace=moment}, for every integer $d$, $\Tr_d(\mathcal{A}(C_l^{(r)}))$ is rational, while $\lambda_{j,k}^2$ may be irrational for $1\leq j\leq l$ and $1\leq k\leq j$. Thus fix any $j$ ($j\in[l]$), $m(j,k_1)=m(j,k_2)$, for any $k_1,k_2\in [j]$. Let $m_j=m(j,1)=m(j,2)=\cdots=m(j,j)$, then the lemma follows directly.$\hfill$$\square$

\noindent\begin{remark}
When $r=3$, the multiplicity $m_l$ does not exist. Thus
\begin{equation}\label{eq:1'}
    \phi(C_{l}^{(3)};\lambda)=\lambda^{m_0}(\lambda^3-\lambda_{1,1}^2)^{m_1}\prod_{j=2}^{l-1}[\tilde{\phi}(P_j;\lambda)]^{m_j}.
\end{equation}
\end{remark}

Observing Equation (\ref{eq:1}), according to Theorem \ref{trace=moment} we have
\begin{equation}\label{eq:mianeqsystem}
    \Tr_{dr}(\mathcal{A}(C_l^{(r)}))=r\sum_{j=1}^{l}m_j\big(\sum_{k=1}^{j}\lambda_{j,k}^{2d}\big),
\end{equation}
which is a system of linear equations in terms of $m_j$ ($j\in[l]$). Then by giving $l$ higher-order traces $\Tr_{dr}(\mathcal{A}(C_l^{(r)}))$ ($d\in[l]$), we can solve for the unknown parameters $m_j$ ($j\in[l]$), and thus obtain $ \phi(C_{l}^{(r)};\lambda)$.

\section{Higher-order traces of $\mathcal{A}(C_{l}^{(r)})$}
\label{sec:ch-inco}
In this section, we derive higher-order traces $\Tr_{dr}(\mathcal{A}(C_l^{(r)}))$ for $d=1,2,\ldots,l$. BEST Theorem \cite{AB} and Matrix-Tree Theorem \cite{CvRS} play an important role in this process since it involves calcutaling numbers of Eulerian cycles and spanning trees of a multi-digraph $D$. For a multi-digraph $D$, denote by $\mathfrak{E}(D)$ (resp. $\tau(D)$) its set of Eulerian circuits (resp. its number of spanning trees). The Laplacian Matrix of $D$ is defined by $L(D) = Deg(D) - A(D)$, where
$Deg(D)$ is a diagonal matrix whose entries are the out-degrees of the vertices of $D$.

\noindent\begin{theorem}\label{th:best}(\cite{AB})
Let $D=(V(D),E(D))$ be an $n$-vertex (multi-)digraph, then
\begin{equation*}
    |\mathfrak{E}(D)|=\tau(D)\prod_{v\in V(D)}\big(d^+(v)-1\big)!.
\end{equation*}
\end{theorem}

\noindent\begin{theorem}\label{th:matrixtree}(\cite{CvRS})
Let $D=(V(D),E(D))$ be an $n$-vertex (multi-)digraph, $L(D)=L$ be its Laplacian Matrix, and $Det(L_{ii})$ be the $i$-th primary minor of $Det(L)$. Then, for any $i\in[n]$
\begin{equation*}
    \tau(D)=Det(L_{ii}).
\end{equation*}
\end{theorem}

It should be noted that in BEST Theorem, multi-arcs are distinct from each other, while in Definition \ref{de:trace2}, they are identified. Thus $|\mathfrak{E}(F)||E(F)|=|W(F)|b(F)$ for any $F\in\mathcal{F}_{dr}'$. Let $p(F)=\prod_{v\in V(F)}d^+(v)$. Further, we have
\begin{equation}\label{eq:newtrace}
\begin{split}
    \Tr_{dr}(\mathcal{T})& = (r-1)^{n-1}\sum_{F\in \mathcal{F}_{dr}'}\frac{b(F)}{c(F)}|W(F)|\pi_{F}(\mathcal{T}) \\
    & = (r-1)^{n-1}\sum_{F\in \mathcal{F}_{dr}'}\frac{E(F)}{c(F)}|\mathfrak{E}(F)|\pi_{F}(\mathcal{T}) \\
    & = (r-1)^{n-1}\sum_{F\in \mathcal{F}_{dr}'}\frac{E(F)}{p(F)}\tau(F)\pi_{F}(\mathcal{T}) \\
     & = (r-1)^{n}\sum_{F\in \mathcal{F}_{dr}'}\frac{dr\tau(F)}{p(F)}\pi_{F}(\mathcal{T}),
\end{split}
\end{equation}
the last equation holds for $|E(F)|=dr(r-1)$.

Due to (3) of Lemma \ref{le:r-symmetry}, we consider $\Tr_{dr}(\mathcal{A}(C_{l}^{(r)}))$ for $d\in[l]$. Let $F=(i_1\alpha_1,\ldots,i_{dr}\alpha_{dr})\in \mathcal{F}_{dr}'$, if $\pi(F)\neq 0$, then every entry $i_s\alpha_s$ of $F$ must form an edge of $C_l^{(r)}$. And if $W(F)\neq 0$, then the multi-digraph $D(F)$ must be connected, consequently the edges involved in $F$ must construct a connected sub-hypergraph in $C_l^{(r)}$, i.e., $P_{s+1}^{(r)}$ ($1\leq s< l$) or $C_l^{(r)}$ itself. For any positive integer $m>1$, let $J_m$ (resp. $J_m^\top$) be an $m$ dimensional all-ones column (resp. raw) vector, and $A_m$ be an $m\times m$ matrix with the following form
\renewcommand{\arraystretch}{1.5}
\begin{equation*}
    A_m=\begin{pmatrix}
    r-1 & -1 & \cdots & -1\\
    -1 & r-1 & \cdots & -1\\
    \vdots & \vdots & \ddots & \vdots\\
    -1 & -1 & \cdots &r-1
\end{pmatrix}.
\end{equation*}

When $d<l$, then hyperedges involved in $F$ construct only hyperpaths with length $s\in [d]$. In this case, we present the primary minor $Det(L_{11})$ (denoted by $p_s$) of $Det(L(D))$, where $L(D)$ is the Laplacian Matrix of multi-digraph $D=D(F)$.
\setlength{\arraycolsep}{1pt}
\begin{equation}\label{eq:p_s}
    p_s=\begin{vmatrix}
a_1A_{r-2} & -a_1J_{r-2} \\
-a_1J_{r-2}^\top & x_1 & -a_2J_{r-2}^\top & -a_2 \\
 & -a_2J_{r-2} & a_2A_{r-2} & -a_2J_{r-2} \\
 &-a_2 & -a_2J_{r-2}^\top & x_2&  \\
 & & & & \ddots \\
 & & & & & x_{s-2} & -a_{s-1}J_{r-2}^\top & -a_{s-1} \\
 & & & & &-a_{s-1}J_{r-2} & a_{s-1}A_{r-2} & -a_{s-1}J_{r-2} &  \\
 & & & & & a_{s-1} & -a_{s-1}J_{r-2}^\top & -x_{s-1} & -a_sJ_{r-1}^\top \\
&&&&&&&-a_sJ_{r-1} & a_sA_{r-1}
\end{vmatrix},
\end{equation}
where $\sum_{i=1}^s a_i=d$, $a_i$ is an nonnegative integer and $x_i=(a_i+a_{i+1})(r-1)$. It should be noted that the total number (including possible multiplicities) of hyperedges involved in $F$ is $dr$, and $a_ir$ refers to the multiplicity of the $i$-th hyperedge of $P_{s+1}^{(r)}$.

When $d\geq l$, then hyperedges involved in $F$ may form $C_{l}^{(r)}$ itself and hyperpaths as well. The case where hyperedges form hyperpaths has been analyzed in the previous context. While in the case where $C_{l}^{(r)}$ is formed, the corresponding multi-digraphs $D(F)$ are not unique; instead, there are two. We present their corresponding primary minors (denoted by $c_l$ and $c'_l$, respectively) as follows
\setlength{\arraycolsep}{1pt}
\begin{equation}\label{eq:c_l}
    c_l=\begin{vmatrix}
a_1A_{r-2} & -a_1J_{r-2} \\
-a_1J_{r-2}^\top & x_1 & -a_2J_{r-2}^\top & -a_2 \\
 & -a_2J_{r-2} & a_2A_{r-2} & -a_2J_{r-2} \\
 &-a_2 & -a_2J_{r-2}^\top & x_2&  \\
 & & & & \ddots \\
 & & & & & x_{l-2} & -a_{l-1}J_{r-2}^\top & -a_{l-1} \\
 & & & & &-a_{l-1}J_{r-2} & a_{l-1}A_{r-2} & -a_{l-1}J_{r-2} &  \\
 & & & & & a_{l-1} & -a_{l-1}J_{r-2}^\top & -x_{l-1} & -a_lJ_{r-2}^\top \\
&&&&&&&-a_lJ_{r-2} & a_lA_{r-2}
\end{vmatrix},
\end{equation}
where $a_1=\cdots =a_l=1$, when $d=l$.
\setlength{\arraycolsep}{1.5pt}
\begin{equation}\label{eq:c'_l}
c'_l=a^{(r-1)l-1}\begin{vmatrix}
 A_{r-2} & -J_{r-2} \\
 & y & -2J_{r-2}^\top & -2 \\
 &  -J_{r-2} & A_{r-2} & -J_{r-2} \\
 &&&y \\
 & & & & \ddots \\
 & & & & & y & -2J_{r-2}^\top \\
 & & & & &  -J_{r-2} & A_{r-2}
\end{vmatrix},\end{equation}
where $y=2(r-1)$ and $la=d$. Note that $ar$ is the common multiplicity of every edge. The order of $c'_l$ is $(r-1)l-1$.

\noindent\begin{lemma}\label{le:plcl}
    Let $p_s,$ $c_l$ and $c'_l$ be determinants shown in Equations (\ref{eq:p_s}--\ref{eq:c'_l}), respectively, and $l,r\geq3$. Then
\begin{align}
    p_s & = r^{s(r-2)}\prod_{i=1}^{s}a_i^{r-1}, \label{eq:p_svalue}\\
    c_l & = 2r^{l(r-2)-1}\prod_{i=1}^{l}a_i^{r-1}\big(\sum_{i=1}^{l}\frac{1}{a_i}\big), \label{eq:c_lvalue}\\
    c'_l & = 2^{l}r^{l(r-2)-1}a^{(r-1)l-1}. \label{eq:c'_lvalue}
\end{align}
\end{lemma}

\noindent\textbf{Proof.} These results can be obtained through two rounds of induction, respectively. First, set $r=3$ and perform induction on $l$, and then perform induction on $r$. $\hfill$$\square$

\noindent\begin{theorem}\label{th:3-1}
Let $C_{l}^{(r)}$ be an $r(\geq 3)$-uniform hypercycle with length $l$, $\mathcal{A}(C_{l}^{(r)})$ be its adjacency tensor. For any positive integer $d$ and any positive integer $s\leq d$, Let
\begin{equation*}
    h(d;s)=\begin{cases}
        1, \quad & s=1, \\
        d\bigg[\sum_{a_1,a_2,\ldots,a_s\geq1}^{a_1+a_2+\cdots+a_s=d}\frac{\prod_{j=2}^{s-1}a_j\prod_{j=1}^{s-1}\binom{a_j+a_{j+1}}{a_j}}{\prod_{j=1}^{s-1}(a_j+a_{j+1})}\bigg], & 2\leq s\leq d.
    \end{cases}
\end{equation*}
Then for every integer $d$ ($1\leq d<l$),
\begin{equation*}
   \Tr_{dr}(\mathcal{A}(C_{l}^{(r)})) = \sum_{s=1}^{d}h(d;s)lr^{s(r-2)+1}(r-1)^{(l-s)(r-1)-1}.
\end{equation*}
\end{theorem}

\noindent\textbf{Proof.}
By Definition \ref{de:trace2}, take $F=(i_{1}\alpha_{1},i_{2}\alpha_2 \ldots, i_{dr}\alpha_{dr})\in \mathcal{F}'_{dr}$. As we have elaborated above, every entry $i\alpha$ in $F$ corresponds to a hyperedge of $C_l^{(r)}$. Suppose that there are $s$ distinct hyperedges $e_j=\{v_1^{(j)},v_2^{(j)},\ldots,v_r^{(j)}\}$ ($j\in[s]$) involved. We map every $i\alpha$ to $ve$, where $i=v\in V(C_{l}^{(r)})$ and $e$ containing $i=v$ is the hyperedge corresponding to $i\alpha$. And we abbreviate an $a$-tuple $ii\cdots i\in [n]^a$ as $i^a$. Since $|W(F)|\neq0$, $F$ is $r$-valent. Specifically, if $v\in V(C_{l}^{(r)})$ appears $a$ times as the first component of some $i\alpha$, then $v$ appears $a(r-1)$ times in these $\alpha$'s. Thus $F$ maps to $\big((v_1^{(1)}e_1)^{a_1}(v_2^{(1)}e_1)^{a_1}\cdots (v_r^{(1)}e_1)^{a_1}(v_1^{(2)}e_2)^{a_2}\cdots (v_r^{(s)}e_s)^{a_s}\big)$, where $v_1^{(1)}<v_2^{(1)}<\ldots<v_r^{(j)}=v_1^{(j+1)}<\ldots<v_r^{(s)}$ and $\sum_{j=1}^s a_j=d$. For any $j\in[s-1]$, $v_r^{(j)}=v_1^{(j+1)}$ is due to that $D(F)$ is connected. The multi-digraph $D(F)$ is indirectly demonstrated by the primary minor $p_s=Det(L_{11})$ in Equation (\ref{eq:p_s}).

When $s=1$, there are $l$ choices for the one hyperedge involved in $F$. Fix a hyperedge and a vertex as the first index of the corresponding $i\alpha$, there are $(r-1)!$ permutations of vertices to form the $(r-1)$-tuple $\alpha$. Thus the total number of such $F$'s is $l[(r-1)!]^{dr}=\frac{l}{\pi_{F}(\mathcal{A}(C_l^{(r)}))}$. And $p(F)=[d(r-1)]^r$. Thus the contribution of any such $F$ to $\Tr_{dr}(\mathcal{A}(C_{l}^{(r)}))$ is
\begin{equation*}
    lr^{r-1}(r-1)^{r(l-1)-l}.
\end{equation*}

When $2\leq s\leq d$, the $s$ hyperedges construct a hyperpath (a sub-hypergraph of $C_l^{(r)}$) with length $s$, there are $l$ such sub-hypergraphs in $C_l^{(r)}$. Fix a hyperedge and a vertex as the first index of the corresponding $i\alpha$, there are $(r-1)!$ permutations of vertices to form $(r-1)$-tuple $\alpha$. Since $v_r^{(j)}=v_1^{(j+1)}$, there are $\binom{a_j+a_{j+1}}{a_j}$ permutations of $v_r^{(j)}e_j$ and $v_1^{(j+1)}e_{j+1}$, for $1\leq j<s$. Thus the total number of $F\in\mathcal{F}'_{dr}$ which involves such $s$ hyperedges is
\begin{equation*}
    l[(r-1)!]^{dr}\bigg[\sum_{a_1,\ldots,a_s\geq1}^{a_1+\cdots+a_s=d}\prod_{j=1}^{s-1}\binom{a_k+a_{k+1}}{a_k}\bigg]=\frac{l}{\pi_F(\mathcal{A}(C_l^{(r)}))}\bigg[\sum_{a_1,\ldots,a_s\geq1}^{a_1+\cdots+a_s=d}\prod_{j=1}^{s-1}\binom{a_k+a_{k+1}}{a_k}\bigg].
\end{equation*}

According to the structure of $D(F)$,
\begin{equation*}
\begin{split}
    p(F) =& \prod_{v\in V(F)}d^+(v) \\
    =& [a_1(r-1)]^{r-1}[a_2(r-1)]^{r-2}\cdots[a_{s-1}(r-1)]^{r-2}[a_s(r-1)]^{r-1}\\
    &\cdot(a_1+a_2)\cdots (a_{s-1}+a_s)(r-1)^{s-1} \\
    =& a_1a_s\bigg[\prod_{j=1}^s a_j^{r-2}\bigg]\bigg[\prod_{j=1}^{s-1}(a_j+a_{j+1})\bigg] (r-1)^{s(r-1)+1}.
\end{split}
\end{equation*}
Then by Equations (\ref{eq:newtrace}--\ref{eq:p_s}, \ref{eq:p_svalue}), the contribution of any such $F$ to $\Tr_{dr}(\mathcal{A}(C_{l}^{(r)}))$ is
\begin{equation*}
dlr^{s(r-2)+1}(r-1)^{(l-s)(r-1)-1}\bigg[\sum_{a_1,a_2,\ldots,a_s\geq1}^{a_1+a_2+\cdots+a_s=d}\frac{\prod_{j=2}^{s-1}a_j\prod_{j=1}^{s-1}\binom{a_j+a_{j+1}}{a_j}}{\prod_{j=1}^{s-1}(a_j+a_{j+1})}\bigg].
\end{equation*}
Note that $\prod_{j=2}^{s-1}a_j=1$, when $s=2$.

Since $1\leq s\leq d$, then we complete the proof $\hfill$$\square$.

\noindent\begin{theorem}\label{th:3-2}
Let $C_{l}^{(r)}$ be an $r(\geq 3)$-uniform hypercycle with length $l$, $\mathcal{A}(C_{l}^{(r)})$ be its adjacency tensor, then
\begin{equation*}
\Tr_{lr}(\mathcal{A}(C_{l}^{(r)})) = \sum_{s=1}^{l-1}h(l;s)lr^{s(r-2)+1}(r-1)^{(l-s)(r-1)-1}+2(l+1)lr^{l(r-2)}.
\end{equation*}
\end{theorem}
\noindent\textbf{Proof.}
Consider $F=(i_{1}\alpha_{1},i_2\alpha_2, \ldots, i_{lr}\alpha_{lr})\in \mathcal{F}'_{lr}$. Similarly, suppose that all the $i\alpha$'s in $F$ correspond totally to $s$ distinct hyperedges constructing a connected sub-hypergraph of $C_{l}^{(r)}$.

When $1\leq s<l$, such a sub-hypergraph can only be a hyperpath of length $s$. And this case has already been discussed in Theorem \ref{th:3-1}, then the total contribution of such $F$'s is $\sum_{s=1}^{l-1}h(l;s)$.

When $s=l$, such a sub-hypergraph can only be $C_l^{(r)}$ itself. To ensure that $|W(F)|\neq0$, every vertex $v\in V(F)$ must have equal in-degree and out-degree, i.e., $d^+(v)=d^-(v)$. Under the constraint that $|W(F)|\neq0$, the multi-digraph $D(F)$ must fall into one of the following cases

{\bf Case 1.} $F$ maps to $\big(v_1^{(1)}e_1v_2^{(1)}e_1\cdots v_r^{(1)}e_1v_1^{(2)}e_2\cdots v_r^{(l)}e_l\big)$, where $v_r^{(l)}=v_1^{(1)}<v_2^{(1)}<\cdots<v_r^{(j)}=v_1^{(j+1)}<\cdots<v_{r-1}^{(l)}$. And the structure of the corresponding multi-digraph $D(F)$ is indirectly shown by its primary minor $c_l=Det(L_{11})$ in Equation \ref{eq:c_l}, where $a_1=\cdots=a_l=1$.

There are $(r-1)!$ permutations of vertices to form each $\alpha$ in $F$. Since $v_r^{(j)}=v_1^{(j+1)}$ ($j\in[l-1]$) and $v_r^{(l)}=v_1^{(1)}$. For each pair (i.e., $(v_r^{(j)},v_1^{(j+1)})$ or $(v_1^{(1)},v_r^{(l)})$) that appears in $F$, there are two possible arrangements (i.e., $v_r^{(j)}e_jv_1^{(j+1)}e_{j+1}$ and $v_r^{(j)}e_{j+1}v_1^{(j+1)}e_{j}$) to choose from. Thus the number of such $F$'s is
\begin{equation*}
    2^{l}[(r-1)!]^{lr}=\frac{2^l}{\pi_{F}(\mathcal{A}(C_{l}^{(r)}))}.
\end{equation*}
By Equation (\ref{eq:c_lvalue}), for every such $F$, $\tau(F)=2lr^{l(r-2)-1}$, and $p(F)=(r-1)^{l(r-2)}(2r-2)^{l}$.
Therefore the total contribution of these $F$'s to $\Tr_{lr}(\mathcal{A}(C_{l}^{(r)}))$ is
\begin{equation*}
    (r-1)^{l(r-1)}\frac{lr\cdot 2lr^{l(r-2)-1}\cdot 2^l}{(r-1)^{l(r-2)}(2r-2)^{l}}=2l^2r^{l(r-2)}.
\end{equation*}

{\bf Case 2.} $F$ maps to $\big((v_1^{(1)}e_1)^{2}v_2^{(1)}e_1\cdots v_{r-1}^{(1)}e_1(v_1^{(2)}e_2)^{2}\cdots (v_1^{(l)}e_l)^2\cdots v_{r-1}^{(l)}e_l\big)$, where $v_r^{(l)}=v_1^{(1)}<v_2^{(1)}<\cdots<v_r^{(j)}=v_1^{(j+1)}<\cdots<v_{r-1}^{(l)}$. And the structure of the corresponding multi-digraph $D(F)$ is indirectly shown by its primary minor $c'_l=Det(L_{11})$ in Equation (\ref{eq:c'_l}), where $a=1$.

There are $(r-1)!$ permutations of vertices to form each $\alpha$ in $F$. Since $v_r^{(j)}=v_1^{(j+1)}\in e_j,e_{j+1}$ ($j\in[l-1]$) and $v_r^{(l)}=v_1^{(1)}\in e_1,e_l$, there are two possible arrangements (i.e., $(v_r^{(j)}e_j)^2$ and $(v_r^{(j)}e_{j+1})^2$) for every such vertex to choose from. Thus the number of such $F$'s is $2[(r-1)!]^{lr}$. By Equation (\ref{eq:c'_lvalue}), $\tau(F)=2^lr^{l(r-2)-1}$, and $p(F)=(r-1)^{l(r-2)}(2r-2)^{l}$. Therefore the total contribution of these $F$'s to $\Tr_{lr}(\mathcal{A}(C_{l}^{(r)}))$ is
\begin{equation*}
    (r-1)^{l(r-1)}\frac{lr\cdot 2^lr^{l(r-2)-1}\cdot 2}{(r-1)^{l(r-2)}(2r-2)^{l}}=2lr^{l(r-2)}.
\end{equation*}

Therefore, summing up all the contributions yields
\begin{equation*}
\Tr_{lr}(\mathcal{A}(C_{l}^{(r)})) = \sum_{s=1}^{l-1}h(l;s)lr^{s(r-2)+1}(r-1)^{(l-s)(r-1)-1}+2(l+1)lr^{l(r-2)}.
\end{equation*}$\hfill\square$

\noindent\begin{corollary}\label{co:tracetotal}
Let $\boldsymbol{T}=\big(\Tr_{r}(\mathcal{A}(C_{l}^{(r)})),\Tr_{2r}(\mathcal{A}(C_{l}^{(r)})),\ldots,\Tr_{lr}(\mathcal{A}(C_{l}^{(r)}))\big)^{\top}$, $\boldsymbol{t}=\big(lr^{r-1}(r-1)^{(l-1)(r-1)-1},lr^{2r-3}(r-1)^{(l-2)(r-1)-1},\ldots,lr^{(l-1)(r-2)+1}(r-1)^{r-2},lr^{l(r-2)}\big)^{\top}$ and $H=(h_{ij})$ be an $l\times l$ lower triangular matrix with the entry
\begin{equation*}
    h_{ij}=\begin{cases}
         0, \quad & j>i, \\
         2(l+1), & j=i=l, \\
        h(i;j), & \text{otherwise}.
    \end{cases}
\end{equation*}
Then combining with Theorems \ref{th:3-1} and \ref{th:3-2} we have
\begin{equation*}
    H\boldsymbol{t}=\boldsymbol{T}.
\end{equation*}
\end{corollary}

\section{The characteristic polynomial of $C_l^{(r)}$}
\label{sec:spectralmoment}

In this section, we solve for the unknown parameters $m_0,m_1,\ldots,m_l$ in Equation (\ref{eq:1}). Let $C_l^{(r)}$ be an $r(\geq 4)$-uniform hypercycle with length $l$. $\lambda_{1,1}=2$ and $\lambda_{j,k}=2\cos{\frac{k\pi}{j+1}}$ is the $k$-th eigenvalue of the $j$-vertex path $P_j$ for $1\leq k\leq j$ and $2\leq j\leq l$. First, we reformulate the system of linear equations (\ref{eq:mianeqsystem}) into matrix form
\begin{equation}\label{eq:matrixform}
    \boldsymbol{T}=rS\boldsymbol{m},
\end{equation}
where $\boldsymbol{T}$ is defined in Corollary \ref{co:tracetotal}, $\boldsymbol{m}=(m_1,m_2,\ldots,m_l)^{\top}$ and $S=(S_{ij})$, $S_{ij}=\sum_{k=1}^j \lambda_{j,k}^{2i}$.

Let $\sigma_k(x_1,x_2,\ldots,x_l)=\sigma_k$ be the $k$-th elementary symmetric polynomial in the polynomial ring $K[x_1,x_2,\ldots,x_l]$ ($0\leq k \leq l$), i.e., $\sigma_0=1$ and $\sigma_{k}=\sum\limits_{1\leq i_1\leq i_2\leq\ldots\leq i_k\leq l}x_{i_1}x_{i_2}\cdots x_{i_k}$ ($1\leq k\leq l$).
Specially, define $\sigma_k(\bar{x}_j)=\sigma_k(x_1,x_2,\ldots,x_{j-1},x_{j+1},\ldots,x_l)$ for $j\in[l]$. Let $V$ be an $l\times l$ Vandermonde matrix, i.e.,
\begin{equation*}
    U=\begin{pmatrix}
    1 & 1 & \cdots & 1\\
    x_1 & x_2 & \cdots & x_l\\
    \vdots & \vdots & \ddots & \vdots\\
    x_1^{l-1} & x_2^{l-1} & \cdots & x_l^{l-1}
\end{pmatrix}.
\end{equation*}

\noindent\begin{lemma}\label{le:Vandermonde}
Let $U^{-1}=(w_{ij})$ be the inverse of $U$. Then
\begin{equation*}
\begin{array}{ccc}
    w_{ij} & = & \frac{(-1)^{l-i}\sigma_{l-i}(\bar{x}_j)}{\prod_{k=1(k\neq j)}^{l}(x_j-x_k)}, \\
    Det(U)&= & \prod_{1\leq i<j\leq l}(x_i-x_j).
\end{array}
\end{equation*}
\end{lemma}

The inverse of the matrix $S$ is given as follows.

\noindent\begin{theorem}\label{th:inverse}
Let $\lambda_{1,1}=2$, $\lambda_{j,k}=2\cos{\frac{k\pi}{j+1}}$ and $S=(S_{ij})=(\sum_{k=1}^j \lambda_{j,k}^{2i})$, $S^{-1}=(S^{-1}_{ij})$ be the inverse of $S$. Let $\alpha\in [1]\times[2]\times\cdots\times[l]$ be an $l$-tuple, i.e., $\alpha_i\in[i]$ for $1\leq i \leq l$. Then
\begin{equation*}
S^{-1}_{ij}=(-1)^{i+j-l+1}\sum_{\alpha\in[1]\times\cdots\times[l]}\sigma_{l-1}(\bar{\lambda}_{i,\alpha_i}^2)\sigma_{l-i}(\bar{\lambda}_{j,\alpha_j}^2)\prod_{\substack{1\leq k<h\leq l\\ k,h\neq j}}(\lambda_{h,\alpha_h}^2-\lambda_{k,\alpha_k}^2)/[2^{l}(l+1)!],
\end{equation*}
where $\sigma_{l-i}(\bar{\lambda}_{j,\alpha_j}^2)=\sigma_{l-i}(\lambda_{1,\alpha_1}^2,\ldots,\lambda_{j-1,\alpha_{j-1}}^2,\lambda_{j+1,\alpha_{j+1}}^2,\ldots,\lambda_{l,\alpha_{l}}^2)$.
\end{theorem}

\noindent\textbf{Proof.} Fix an $\alpha\in [1]\times[2]\times\cdots\times[l]$, let
\setlength{\arraycolsep}{3.pt}
\begin{equation*}
M=\begin{pmatrix}
    \lambda_{1,\alpha_1}^2 & \lambda_{2,\alpha_2}^2 & \cdots & \lambda_{l,\alpha_l}^2 \\
    \lambda_{1,\alpha_1}^4 & \lambda_{2,\alpha_2}^4 & \cdots & \lambda_{l,\alpha_l}^4 \\
    \vdots & \vdots & \ddots & \vdots \\
    \lambda_{1,\alpha_1}^{2l} & \lambda_{2,\alpha_2}^{2l} & \cdots & \lambda_{l,\alpha_l}^{2l}
\end{pmatrix}.
\end{equation*} According to Lemma \ref{le:Vandermonde} and the multilinearity of the determinant in columns, we have
\begin{equation*}
    Det(M)=\sigma_l(\lambda_{1,\alpha_1}^2,\ldots,\lambda_{l,\alpha_l}^2)\prod_{1\leq k<h\leq l}(\lambda_{h,\alpha_h}^2-\lambda_{k,\alpha_k}^2)=\sigma_l\prod_{1\leq k<h\leq l}(\lambda_{h,\alpha_h}^2-\lambda_{k,\alpha_k}^2).
\end{equation*}
Furthermore, the multilinearity of the determinant in the columns of $S$ yields
\begin{equation*}
    Det(S)=\sum_{\alpha\in[1]\times\cdots\times[l]}\sigma_l\prod_{1\leq k<h\leq l}(\lambda_{h,\alpha_h}^2-\lambda_{k,\alpha_k}^2),
\end{equation*}
where $\sigma_l=\sigma_l(\lambda_{1,\alpha_1}^2,\ldots,\lambda_{l,\alpha_l}^2)=\prod_{i=1}^{l}\lambda_{i,\alpha_i}^2$. By induction on $l$ combined with computer-aided computation, we have
\begin{equation*}
    Det(S)=(-1)^{l-1}2^{l}(l+1)!.
\end{equation*}

Denote by $M^*=(M^*_{ij})$ the adjugate matrix of $M$. For the fixed
$l$-tuple $\alpha$, identify $x_t$ in $U$ with $\lambda_{t,\alpha_t}$ for each $t\in[l]$. Then
\begin{equation*}
\begin{split}
M^*_{ij} &= \prod_{t=1(t\neq i)}^{l}x_t U^*_{ij} \\
&= \sigma_{l-1}(\bar{x}_{i})(w_{ij}Det(U_l)) \\
&= \sigma_{l-1}(\bar{x}_{i})\prod_{1\leq k<h\leq l}(x_h-x_k)\frac{(-1)^{l-i}\sigma_{l-i}(\bar{x}_j)}{\prod_{k=1(k\neq j)}^{l}(x_j-x_k)} \\
&= (-1)^{n-i+n-j}\prod_{\substack{1\leq k<h\leq l\\k,h\neq j}}(x_h-x_k)\sigma_{l-i}(\bar{x}_j)\sigma_{l-1}(\bar{x}_i)\\
&= (-1)^{i+j}\sigma_{l-i}(\bar{x}_j)\sigma_{l-1}(\bar{x}_i)\prod_{\substack{1\leq k<h\leq l\\k,h\neq j}}(x_h-x_k).
\end{split}
\end{equation*}

Similarly, by the multilinearity of the determinant in columns of $S$ we have
\begin{equation*}
    S^{*}_{ij}=(-1)^{i+j}\sum_{\alpha\in[1]\times\cdots\times[l]}\sigma_{l-1}(\bar{\lambda}_{i,\alpha_i}^2)\sigma_{l-i}(\bar{\lambda}_{j,\alpha_j}^2)\prod_{\substack{1\leq k<h\leq l\\ k,h\neq j}}(\lambda_{h,\alpha_h}^2-\lambda_{k,\alpha_k}^2).
\end{equation*}
Since $S^{-1}_{ij}=S^{*}_{ij}/Det(S)$, then Theorem \ref{th:inverse} follows.$\hfill$$\square$

\noindent\begin{theorem}\label{th:polynomial}
Let $C_{l}^{(r)}$ be an $r(\geq 4)$-uniform hypercycle with length $l$. Then
\begin{equation*}
        \phi(C_{l}^{(r)};\lambda)=\lambda^{m_0}(\lambda^r-4)^{m_1}\prod_{j=2}^{l}\tilde{\phi}(P_j;\lambda)^{m_j},
\end{equation*}
where
\begin{equation*}
\begin{split}
    m_i &= r^{-1}\sum_{j=1}^{l}S^{-1}_{ij}\Tr_{jr}(\mathcal{A}(C_{l}^{(r)})) \quad(i\in[l]), \\
    m_0 &= l(r-1)^{l(r-1)}-\sum_{i=1}^{l}irm_i.
\end{split}
\end{equation*}
and $\Tr_{jr}(\mathcal{A}(C_{l}^{(r)}))$, $S^{-1}_{ij}$ are shown in Theorems \ref{th:3-1}, \ref{th:3-2} and \ref{th:inverse}, respectively.
\end{theorem}

\noindent\textbf{Proof.} By Equation (\ref{eq:matrixform}), $m_i=r^{-1}\sum_{j=1}^{l}S^{-1}_{ij}\Tr_{jr}(\mathcal{A}(C_{l}^{(r)}))$ for $1\leq i\leq l$. As introduced in Definition \ref{de:defofcharpoly}, the total degree of $\phi(C_{l}^{(r)};\lambda)$ is $n(r-1)^{n-1}$, where $n=l(r-1)$. Thus, $m_0=l(r-1)^{l(r-1)}-\sum_{i=1}^{l}irm_i$. From Equation (\ref{eq:1}), the validity of the theorem follows directly.$\hfill$$\square$

The method of solving for $\boldsymbol{m}$ by finding the inverse of matrix $S$ as described in Theorem \ref{th:inverse} is relatively complex. Next, we present a more straightforward and concise approach to determine $\boldsymbol{m}$. Let $B$ be an $l\times l$ matrix as follows
\begin{equation}\label{eq:B}
    B=\begin{pmatrix}
 &  &  &  &  &  & \frac{1}{l}\\
\frac{1}{2} & -1 & \frac{1}{2} &  &  &  & \\
 & \frac{1}{2} & -1 & \frac{1}{2} &  &  & \\
 &  & \ddots  & \ddots  & \ddots  &  & \\
 &  &  & \frac{1}{2} & -1 & \frac{1}{2} & \\
 &  &  &  & \frac{1}{2} & -1 & \frac{l+1}{l}\\
 &  &  &  &  & \frac{1}{2} & -1
\end{pmatrix}.
\end{equation}
It is easy to check that
\begin{equation}\label{eq:inverseB}
    B^{-1}=\begin{pmatrix}
4 & 2 & 4 & 6 & 8 & \cdots  & 2( l-1)\\
6 &  & 2 & 4 & \ddots  & \ddots  & \vdots \\
8 &  &  & 2 & \ddots  & 6 & 8\\
\vdots  &  &  &  & \ddots  & 4 & 6\\
2( l-1) &  &  &  &  & 2 & 4\\
2l &  &  &  &  &  & 2\\
l &  &  &  &  &  &
\end{pmatrix}.
\end{equation}

\noindent\begin{theorem}\label{th:S=HB^-1}
Let $S$ (resp. $B^{-1}$) be defined in Equation (\ref{eq:matrixform}) (resp. Equation (\ref{eq:inverseB})) and $H$ be defined in Corollary \ref{co:tracetotal}. Then $S=HB^{-1}$, i.e.,
\begin{equation*}
    S_{ij}=\begin{cases}
        \sum_{k=1}^{j}\lambda_{j,k}^{2i}=\sum_{s=1}^{\min\{i,j-1\}}2(j-s)h(i;s), \quad & 2\leq j\leq l, \\
        4^{i}=\sum_{s=1}^{i}2(s+1)h_{is}, & j=1,1\leq i\leq l.
    \end{cases}
\end{equation*}
\end{theorem}

\noindent\textbf{Proof.}Depending on whether $2\leq j$, we consider the following two cases.

\noindent\textbf{Case 1.} $2\leq j\leq l$.

In this case, $S_{ij}=\sum_{k=1}^{j}\lambda_{j,k}^{2i}$, which also equals to the number of closed traces with length $2i$ in the $j$-vertex path $P_j=(V(P_j),E(P_j))$. To streamline the subsequent discussion, let $V(P_j)=\{v_1,v_2\ldots,v_j\}$, $E(P_j)=\{e_1,e_2,\ldots,e_{j-1}\}$ and $e_s=v_sv_{s+1}$ for $s\in[j-1]$. To indicate the directions of edges within these closed traces, two arcs associated with edge $e_s$ are defined as $arc^+_s=(v_s,v_{s+1})$ (called the positive arc) and $arc^-_{s}=(v_{s+1},v_s)$ (called the negative arc) for $s\in[j-1]$. Note that any closed trace can be identified with an arrangement of arcs; moreover, the multiplicity and order of negative arcs are entirely determined the positive arcs. Hence, we only need to consider the arrangement of positive arcs.

In any closed trace, denote by $s$ the number of distinct edges involved.
In addition, denote by $a_k\geq1$ the multiplicity of the $k$-th distinct positive arc (without loss of generality, let it be $arc^+_{k}$) in the closed trace for $k\in[s]$. Note that $\sum_{k=1}^{s}a_k=i$. First, we determine the number of circular permutations of these $2i$ arcs, where consecutive arcs in each permutation remain adjacent along $P_j$. Then, multiplying this number by $2i$ yields the total number of closed traces. The advantage of this method is that it avoids the need to consider the starting point and direction of the closed trace. When $s=1$, it is trivial because such number of circular permutations of these $2i$ arcs is $1$, equal to $h(i;1)$. When $s\geq2$, first, the number of circular permutations of $a_1$ $arc^+_1$'s and  $a_2$ $arc^+_2$'s is
\begin{equation*}
    \frac{\binom{a_1+a_2}{a_1}}{a_1+a_2}.
\end{equation*}
Starting from the above circular permutation of $arc^+_1$ and $arc^+_2$, we insert $a_3$ $arc^+_{3}$'s into this circular permutation, while preserving the relative order between $arc^+_{1}$'s (resp. $arc^+_{2}$'s) and $arc^+_{2}$'s (resp. $arc^+_{3}$'s). Hence the number of circular permutations of $a_1$ $arc^+_1$'s, $a_2$ $arc^+_2$'s and  $a_3$ $arc^+_3$'s is
\begin{equation*}
    \frac{a_2\binom{a_1+a_2}{a_1}\binom{a_2+a_3}{a_2}}{(a_1+a_2)(a_2+a_3)}.
\end{equation*}
The rationale for the factor $a_{2}$ in the numerator is as follows: To insert $a_{3}$ $arc^{+}_{3}$'s while preserving the relative order of $arc^{+}_{2}$'s and $arc^{+}_{3}$'s, we must fix an initial $arc^+_{2}$ as a reference in the circular permutation. And there are $a_2$ choices for the reference $arc^+_{2}$.

Repeating the insertion procedure yields the total number of circular permutations of the $i$ positive arcs as
\begin{equation*}
    \frac{\prod_{k=2}^{s-1}a_k\prod_{k=1}^{s-1}\binom{a_k+a_{k+1}}{a_k}}{\prod_{k=1}^{s-1}(a_k+a_{k+1})}.
\end{equation*}
Finally, inserting the remaining $i$ negative arcs yields all desired circuit permutations with consecutive arcs adjacent in the path $P_j$.

Note that the path $P_j$ contains $j-s$ subpaths each of length $s$ for $1\leq s\leq j-1$. Thus the number of closed traces, i.e., $\sum_{k=1}^{j}\lambda_{j,k}^{2i}$ is equal to
\begin{equation*}
    2i\sum_{s=1}^{\min\{i,j-1\}}(j-s)\sum_{a_1+a_2+\cdots+a_s=i}\frac{\prod_{k=2}^{s-1}a_k\prod_{k=1}^{s-1}\binom{a_k+a_{k+1}}{a_k}}{\prod_{k=1}^{s-1}(a_k+a_{k+1})}=\sum_{s=1}^{\min\{i,j-1\}}2(j-s)h(i;s).
\end{equation*}

\noindent\textbf{Case 2.} $j=1$. Since the definition of the entry $h_{ll}$ of $H$ is different from other entries, we divide the argument into two subcases.

\noindent\textbf{Subcase 2.1.} $j=1$ and $1\leq i< l$.

\textbf{Claim 1.} For the length-$t$ cycle $C_t$, the number of closed traces with length $2i$ $(i<t)$ in $C_t$ equals to
\begin{equation*}
    2i\sum_{s=1}^{i}t\sum_{a_1+a_2+\cdots+a_s=i}\frac{\prod_{k=2}^{s-1}a_k\prod_{k=1}^{s-1}\binom{a_k+a_{k+1}}{a_k}}{\prod_{k=1}^{s-1}(a_k+a_{k+1})}=\sum_{s=1}^{i}2th(i;s).
\end{equation*}
when $j$ is even and $i<\frac{t}{2}$ or $j$ is odd.

\textbf{Proof of Claim 1.} Under the condition that $j$ is even and $i<\frac{t}{2}$ or $j$ is odd, the edges involved in each closed trace of length $2i$ must form a subpath of $C_t$, and the number of such path with the same length $s$ is $t$. Combining this with the analysis in Case 1, Claim 1 follows directly.

\textbf{Claim 2.} For any integers $t\geq3$ and $i$,
\begin{equation*}
    \sum_{s=1}^{i}2th(i;s)=\sum_{k=1}^{t-1}\lambda_{t-1,k}^{2i}+4^i.
\end{equation*}

\textbf{Proof of Claim 2.} Observe that both sides of the above equation are continuous in $t$. Therefore, it suffices to verify the equality when $t$ is odd. According to Claim 1, when $t$ is odd, the left side of the equation equals to the number of $2i$-length closed traces of $C_t$, i.e., the $2i$-th spectral moment of $C_t$. And $\sum_{k=1}^{t-1}\lambda_{t-1,k}^{2i}$ equals to the $2i$-th spectral moment of $P_{t-1}$. The squared eigenvalues of $C_t$ and $P_{t-1}$ are listed in rows 1 and 2 of Table \ref{tab:2}. Observe that the entries in any Column $k<t$ of Table \ref{tab:2} coincide pairwise. This equivalence confirms the identity stated in Claim 2.
\begin{table}[h]
	\centering
    \scriptsize
	\begin{tabular}{|l|l|l|l|l|l|l|l|l|}
		\hline
		4$\cos^2(\frac{2\pi}{t})$ & 4$\cos^2(\frac{4\pi}{t})$ & $\cdots$ & 4$\cos^2(\frac{(t-1)\pi}{t})$ & 4$\cos^2(\frac{(t+1)\pi}{t})$ & 4$\cos^2(\frac{(t+3)\pi}{t})$ & $\cdots$ & 4$\cos^2(\frac{2(t-1)\pi}{t})$ & 4 \\ \hline
		4$\cos^2(\frac{2\pi}{t})$ & 4$\cos^2(\frac{4\pi}{t})$ & $\cdots$ & 4$\cos^2(\frac{(t-1)\pi}{t})$ & 4$\cos^2(\frac{\pi}{t})$      & 4$\cos^2(\frac{3\pi}{t})$     & $\cdots$ & 4$\cos^2(\frac{(t-2)\pi}{t})$  &   \\ \hline
	\end{tabular}
	\caption{Squared eigenvalues of $C_t$ and $P_{t-1}$}
	\label{tab:2}
\end{table}

We now proceed to verify Subcase 2.1. Take $t=i+2$, then
\begin{equation*}
\begin{split}
    4^i &= \sum_{s=1}^{i}2(i+2)h(i;s)-\sum_{k=1}^{i+1}\lambda_{i+1,k}^{2i}\\
        &= \sum_{s=1}^{i}2(i+2)h(i;s)-\sum_{s=1}^{i}2(i+1-s)h(i;s) \\
        &= \sum_{s=1}^{i}2(s+1)h(i;s).
\end{split}
\end{equation*}

\noindent\textbf{Subcase 2.2.} $j=1$ and $i=l$.

According to Subcase 2.1,
\begin{equation*}
\begin{split}
    4^l &= \sum_{s=1}^{l-1}2(s+1)h(l;s)+2(l+1)l\bigg[\sum_{a_1,a_2,\ldots,a_l\geq1}^{a_1+a_2+\cdots+a_l=l}\frac{\prod_{k=2}^{l-1}a_k\prod_{k=1}^{l-1}\binom{a_k+a_{k+1}}{a_k}}{\prod_{k=1}^{l-1}(a_k+a_{k+1})}\bigg]\\
    &=\sum_{s=1}^{l-1}2(s+1)h(l;s)+lh_{ll}.
\end{split}
\end{equation*}

Note the relation between $h(i;j)$ and $h_{ij}$, then we complete the whole proof.$\hfill$$\square$

\noindent\begin{theorem}\label{th:main}
Let $C_{l}^{(r)}$ be an $r(\geq 3)$-uniform hypercycle with length $l$. Let $m$, $B$ and $\boldsymbol{t}$  be defined in Equations (\ref{eq:matrixform}-\ref{eq:B}) and Corollary \ref{co:tracetotal}, respectively. Then
\begin{equation*}
    \boldsymbol{m}=r^{-1}B\boldsymbol{t}.
\end{equation*}
and
\begin{equation*}
        \phi(C_{l}^{(r)};\lambda)=\lambda^{m_0}(\lambda^r-4)^{m_1}\prod_{j=2}^{l}\tilde{\phi}(P_j;\lambda)^{m_j},
\end{equation*}
where $m_i$ is the $i$-th component of $\boldsymbol{m}$,  $m_0=l(r-1)^{l(r-1)}-\sum_{i=1}^{l}irm_i$ and $\tilde{\phi}(P_j;\lambda)=\prod_{k=1}^{j}(\lambda^r-\lambda_{j,k}^2)$ is the transformed polynomial of $\phi(P_{j};\lambda)$.
\end{theorem}

\noindent\textbf{Proof.} Lemma \ref{le:4-5} implies the multiplicity  $m_l$ does not occur for $r=3$. Thus we verify this theorem in the following two cases.

\noindent\textbf{Case 1.} $r\geq 4$. According to Corollary \ref{co:tracetotal} and Equation (\ref{eq:matrixform}), we have $H\boldsymbol{t}=rS\boldsymbol{m}$. Since $S=HB^{-1}$ by Theorem \ref{th:S=HB^-1}, we have $\boldsymbol{m}=r^{-1}B\boldsymbol{t}$.

\noindent\textbf{Case 2.} $r=3$. In this case, we have
\begin{equation*}
    m_l=\frac{l}{2}\cdot 3^{l-1+1}\cdot 2-l\cdot 3^l
=0,\end{equation*} which meet the requirement of the institution. And due to $m_l=0$, equations
\begin{equation*}
    \boldsymbol{T}=H\boldsymbol{t}=rS\boldsymbol{m}=rHB^{-1}\boldsymbol{m}
\end{equation*}
still hold in the case that $r=3$.

Thus Theorem \ref{th:main} follows directly.$\hfill$$\square$

One can check that when $l=3$ and $4$, respectively, our results coincide with the known results, $\phi(C_3^{(r)};\lambda)$ (\cite{D2023}) and $\phi(C_4^{(r)};\lambda)$ (\cite{D2024}).

\noindent\begin{corollary} \label{cor:ch1}
Let $r\geq3$. Then
\begin{equation*}
\begin{split}
\phi(C_5^{(r)};\lambda)=&\lambda^{m_0+rm_3+rm_5}(\lambda^r-1)^{2m_2+2m_5}(\lambda^r-2)^{2m_3}(\lambda^r-3)^{2m_5}\\
&(\lambda^r-4)^{m_1}(\lambda^r-\frac{3+\sqrt{5}}{2})^{2m_4}(\lambda^r-\frac{3-\sqrt{5}}{2})^{2m_4},
\end{split}
\end{equation*}

where
\begin{align*}
m_0 =& 5 (r-1)^{5r-5} - 5 (r-1)^{4r-5} r^{r-1} +
 \frac{5}{2} (r-1)^{3r-4} r^{2r-3} - r^{5r-11} + r^{5r-10} \\
m_1 =& r^{5r-11}, \\
m_2 =& \frac{5}{2}(r-1)^{4r-5} r^{r-2} - 5 (r-1)^{3r-4} r^{2r-4} + \frac{5}{2} (r-1)^{2r-3} r^{3r-6},\\
m_3=& \frac{5}{2} (r-1)^{3r-4} r^{2r-4} - 5 (r-1)^{2r-3} r^{3r-6} + \frac{5}{2} (r-1)^{r-2}r^{4r-8},\\
m_4=& \frac{5}{2} (r-1)^{2r-3} r^{3r-6} -
 5 (r-1)^{r-2} r^{4r-8} + 6 r^{5r-11},\\
m_5=& \frac{5}{2} (r-1)^{r-2} r^{4r-8} - 5 r^{5r-11}.
\end{align*}
\end{corollary}

\noindent\textbf{Proof.} By Theorem \ref{th:main}, we have
\begin{equation*}
\begin{split}
\phi(C_5^{(r)};\lambda)=&\lambda^{m_0}(\lambda^r-4)^{m_1}(\lambda^r-1)^{2m_2}[(\lambda^r-2)^2\lambda^r]^{m_3}\\
&[(\lambda^r-\frac{3+\sqrt{5}}{2})(\lambda^r-\frac{3-\sqrt{5}}{2})]^{2m_4}[(\lambda^r-3)^2(\lambda^r-1)^2\lambda^r]^{m_5},
\end{split}
\end{equation*}where
\begin{equation*}
\begin{pmatrix}
    m_1\\
    m_2\\
    m_3\\
    m_4\\
    m_5
\end{pmatrix}=
\begin{pmatrix}
 0 & 0 & 0 & 0 & \frac{1}{5}\\
\frac{1}{2} & -1 & \frac{1}{2} & 0 & 0\\
 0 & \frac{1}{2} & -1 & \frac{1}{2} & 0\\
 0 & 0 & \frac{1}{2} & -1 & \frac{6}{5}\\
 0 & 0 & 0 & \frac{1}{2} & -1
\end{pmatrix}
\left(\begin{array}{c}
    5 (r-1)^{4r-5} r^{r-2} \\
    5 (r-1)^{3r-4} r^{2r-4}\\
    5 (r-1)^{2r-3} r^{3r-6}\\
    5 (r-1)^{r-2} r^{4r-8}\\
    5 r^{5r-11}
\end{array}\right),
\end{equation*}
and $m_0=5(r-1)^{5(r-1)}-\sum_{i=1}^{5}irm_i$. Thus Corollary \ref{cor:ch1} follows directly.$\hfill$$\square$

\noindent\begin{corollary} \label{cor:ch2}
Let $r\geq3$. Then
\begin{equation*}
\begin{split}
\phi(C_6^{(r)};\lambda)=&\lambda^{m_0+rm_3+rm_5}(\lambda^r-1)^{2m_2+2m_5}(\lambda^r-2)^{2m_3}(\lambda^r-3)^{2m_5}\\
&(\lambda^r-4)^{m_1}(\lambda^r-\frac{3+\sqrt{5}}{2})^{2m_4}(\lambda^r-\frac{3-\sqrt{5}}{2})^{2m_4}\\
&(\lambda^r-4\cos^2(\frac{\pi}{7}))^{2m_6}(\lambda^r-4\cos^2(\frac{2\pi}{7}))^{2m_6}(\lambda^r-4\cos^2(\frac{3\pi}{7}))^{2m_6},
\end{split}
\end{equation*}
where
\begin{align*}
m_0 =& 6 (r-1)^{6 r-6}-6 (r-1)^{5 r-6} r^{r-1}+3 (r-1)^{4 r-5} r^{2 r-3}-r^{6 r-13}+r^{6 r-12},\\
m_1 =& r^{6r-13}, \\
m_2 =& 3 (r-1)^{5r-6} r^{r-2} - 6(r-1)^{4r-5}r^{2r-4} + 3(r-1)^{3r-4}r^{3r-6},\\
m_3=& 3 (r-1)^{4r-5} r^{2r-4} - 6(r-1)^{3r-4}r^{3r-6} + 3(r-1)^{2r-3}r^{4r-8},\\
m_4=& 3 (r-1)^{3r-4} r^{3r-6} - 6(r-1)^{2r-3}r^{4r-8} + 3(r-1)^{r-2}r^{5r-12},\\
m_5=& 3 (r-1)^{2r-3} r^{4r-8} - 6(r-1)^{r-2}r^{5r-12} + 7r^{6r-13},\\
m_6=& 3 (r-1)^{r-2} r^{5r-12} - 6 r^{6r-13}.
\end{align*}
\end{corollary}

\noindent\textbf{Proof.} By Theorem \ref{th:main}, we have
\begin{equation*}
\begin{split}
\phi(C_6^{(r)};\lambda)=&\lambda^{m_0}(\lambda^r-4)^{m_1}(\lambda^r-1)^{2m_2}[(\lambda^r-2)^2\lambda^r]^{m_3}\\
&[(\lambda^r-\frac{3+\sqrt{5}}{2})(\lambda^r-\frac{3-\sqrt{5}}{2})]^{2m_4}[(\lambda^r-3)^2(\lambda^r-1)^2\lambda^r]^{m_5}\\
&[(\lambda^r-4\cos^2(\frac{\pi}{7}))(\lambda^r-4\cos^2(\frac{2\pi}{7}))(\lambda^r-4\cos^2(\frac{3\pi}{7}))]^{2m_6},
\end{split}
\end{equation*}where
\begin{equation*}
\begin{pmatrix}
    m_1\\
    m_2\\
    m_3\\
    m_4\\
    m_5\\
    m_6
\end{pmatrix}=
\begin{pmatrix}
 0 & 0 & 0 & 0 & 0 & \frac{1}{6}\\
\frac{1}{2} & -1 & \frac{1}{2} & 0 & 0 & 0\\
 0 & \frac{1}{2} & -1 & \frac{1}{2} & 0 & 0\\
 0 & 0 & \frac{1}{2} & -1 & \frac{1}{2} & 0\\
 0 & 0 & 0 & \frac{1}{2} & -1 & \frac{7}{6}\\
 0 & 0 & 0 & 0 & \frac{1}{2} & -1
\end{pmatrix}
\left(\begin{array}{c}
    6 (r-1)^{5r-6} r^{r-2}\\
    6 (r-1)^{4r-5} r^{2r-4}\\
    6 (r-1)^{3r-4} r^{3r-6}\\
    6 (r-1)^{2r-3} r^{4r-8}\\
    6 (r-1)^{r-2} r^{5r-12}\\
    6 r^{6r-13}
\end{array}\right),
\end{equation*}
and $m_0=6(r-1)^{6(r-1)}-\sum_{i=1}^{6}irm_i$. Thus Corollary \ref{cor:ch2} follows directly. $\hfill$$\square$

\section*{Declaration of competing interest}
The authors declare that they have no conflict of interest.

\section*{Data availability}
The study has no associated data.

\end{document}